%% file: Minimal_sets.tex
\documentclass[12pt]{scrartcl}
\usepackage[left=3cm,right=3cm, top=3cm,bottom=3cm]{geometry}

\newcommand*{\appendixmore}{%
  \renewcommand*{\othersectionlevelsformat}[1]{%
    \ifthenelse{\equal{##1}{section}}{\appendixname~}{}%
    \csname the##1\endcsname\autodot\enskip}
  \renewcommand*{\sectionmarkformat}{%
    \appendixname~\thesection\autodot\enskip}
}

\usepackage{amssymb,amsmath,amscd,amsthm,braket}
\usepackage{datetime, comment}
\usepackage{amsfonts,graphicx,amssymb}
\usepackage[shortlabels]{enumitem}
\usepackage[unicode]{hyperref}
\usepackage{scrpage2}

\usepackage{tikz, wrapfig, graphicx}
\usetikzlibrary{calc,scopes}
\tikzset{>=latex}

\usepackage{mathtools}
\usepackage[mathletters]{ucs}
\let\mathscr=\mathcal
\usepackage[utf8x]{inputenc}
\usepackage[T1]{fontenc}
\usepackage[english]{babel}

\usepackage[backgroundcolor=white, linecolor=black, textwidth=3cm, disable]{todonotes}

\DeclareUnicodeCharacter{8788}{\coloneqq}
\DeclareUnicodeCharacter{8797}{\coloneqq}

\usepackage{aliascnt}
\newtheorem{theorem}{Theorem}[section]
\newcommand{\mynewtheorem}[2]{%
\newaliascnt{#1}{theorem}
\newtheorem{#1}[#1]{#2}
\newtheorem*{#1*}{#2}
\aliascntresetthe{#1}
\expandafter\def\csname #1autorefname\endcsname{#2}
}

\mynewtheorem{proposition}{Proposition}
\mynewtheorem{corollary}{Corollary}
\mynewtheorem{lemma}{Lemma}
\mynewtheorem{claim}{Claim}
\mynewtheorem{question}{Question}
\theoremstyle{remark}
\mynewtheorem{remark}{Remark}
\theoremstyle{definition}
\mynewtheorem{definition}{Definition}

\date{\today\ \currenttime}

\ohead{\headmark}
\ihead{Semigroups of circle diffeomorphisms}
\automark[section]{section}
\ofoot{}

\newcommand{\Exp}{\mathbb E}

\newcommand{\bbZ}{\mathbb{Z}}
\newcommand{\bbR}{\mathbb{R}}
\newcommand{\Sc}{S^1}
\newcommand{\Pbase}{\mathbf P}
\newcommand{\mI}{\mathcal{I}}
\newcommand{\mM}{\mathcal{M}}
\newcommand{\tM}{\widetilde{\mM}}
\newcommand{\nd}{\mathcal{N}}

\newcommand{\rot}{\mathcal{D}}

\newcommand{\eps}{\varepsilon}

\newcommand{\Gr}{\mathcal{S}}

\newcommand{\Homeo}{\mathrm{Homeo}}
\newcommand{\PSL}{\mathrm{PSL}}
\DeclareMathOperator{\supp}{supp}
\DeclareMathOperator{\Id}{Id}

\DeclareMathOperator{\diam}{diam}

\DeclareMathOperator{\Prob}{P} 
\usepackage{authblk}

\title{Classification of generic semigroup actions of circle diffeomorphisms}
\author{
    V. Kleptsyn%
    \thanks{CNRS, Institut de Recherches Math\'ematiques de Rennes, Campus Beaulieu, 35042 Rennes, France.}
    , Yu. Kudryashov%
    \thanks{Department of Mathematics, Cornell University, NY, USA; 
    		 \newline ORCID ID:00000-0003-4286-9276}
    ,  A. Okunev%
    \thanks{Department of Mathematics, Higher School of Economics, Moscow, Russia.}
}

\usepackage{etoolbox}
\makeatletter
\hypersetup{pdftitle=\@title,pdfauthor={Victor Kleptsyn, Yury Kudryashov, Alexey Okunev}}
\hypersetup{pdftitle=\@title}
\expandafter\let\expandafter\unicode@greek@prefix\expandafter=\@firstoftwoα

\newcommand{\unicode@accent}[2]{
    \ifx\unicode@greek@prefix#2
        \let\temp@a=\unicode@accent@two
    \else
        \let\temp@a=\relax
    \fi
    \temp@a#1#2
}
\newcommand{\unicode@accent@two}[3]{#1{#2#3}}
\newcommand{\unicode@patch@accent}[1]{
    \cslet{unicode@save@\string#1}{#1}
    \def#1{\@ifnextchar\bgroup{\@nameuse{unicode@save@\string#1}}{\expandafter\unicode@accent\csname unicode@save@\string#1\endcsname}}
}
\unicode@patch@accent\tilde
\unicode@patch@accent\hat

\makeatother

\begin{document}
\pagenumbering{arabic}

\maketitle

\abstract{We study topological properties of semi-group actions on the circle by orientation-preserving homeomorhisms. We prove that a generic action either possesses a forward-invariant interval-domain (i.e. a finite union of disjoint circle arcs), or is two-sided minimal (and, moreover, by a small perturbation one can create a map with arbitrary rotation number in the semigroup). For the minimal case, we also study conditions required for global synchronization.}

\tableofcontents

\input introduction
\input preliminaries
\input questions

\input topology

\input perturbation

\input minimality

\input factorization
\input synchronization
\input examples

\section*{Acknowledgments}
The authors are grateful to Michele Trienstino, Dominique Malicet and Anton Gorodetski
for many fruitful discussions.
\appendix
\input robust

\input multipliers

\end{document}

%% file: introduction.tex
\section{Introduction}
\subsection{History}
\todo{Cite~\cite{Du}, \cite{BR} in the journal version?}
The purpose of this note is to obtain a topological description for the possible properties of  semi-group actions on the circle by orientation-preserving homeomorhisms, or, what is the same, of the iterated function system (IFS) formed by the generators of the semigroup.
To obtain this description we will use as one of the main tools the random dynamics, choosing some positive probabilities of the generators of the semigroup.

One of the important tools that we use throughout this paper is the (local) contraction. Over the long years of studies of random dynamical systems many authors have remarked that random orbits starting at different initial points approach each other. One of the pioneering works where this phenomemon have been observed was a paper~\cite{F} by H.~Furstenberg, who has established (under very mild assumptions) a random contraction effect for the projectivization of a linear dynamics.

P.~Baxendale's theorem~\cite{B} shows that for a random dynamics by diffeomorphisms on a compact manifold, in absence of a common invariant measure, there exists a stationary measure with negative “volume” Lyapunov exponent. For the one-dimensional dynamics, this implies the local contraction; such an argument was used in~\cite{DKN} for the symmetric random walks for group actions on the circle. In V.~A.~Antonov's work~\cite{A} (see also~\cite{KN}), the smoothness assumptions were removed, at the cost of assuming forward and backward minimality of the action.

Random dynamics on an interval was also studied, turning out to possess even more regular properties (see~\cite{Ku,KV}) than for the case of a circle. In particular, a contraction effect was also noticed to take place. Finally, in a recent outstanding work of D. Malicet~\cite{M}, the effect of local (exponential) contraction in random dynamics on the circle was  established under no assumption other than the absence of a common invariant measure. Among his results is that under these assumptions, almost surely a random trajectory tends to one of the forward-minimal sets~(\cite[Theorem~B]{M}).

\subsection{Main results}
For shortness, we introduce the following definition.
\begin{definition}
    A~union~$U$ of~finitely many open intervals, $U  \notin\set{\Sc,\varnothing}$, is called an~\emph{interval-domain}. An interval-domain $U$ is \emph{(forward-) invariant} for an IFS $(f_1, \dots, f_s)$, if $f_i(U) \subset U$ for each $i$. An interval-domain $U$ is \emph{backward-invariant}, if $f_i^{-1}(U) \subset U$ for each $i$.   
\end{definition}
\noindent If the~action of a~semigroup $G \subset \Homeo(\Sc)$ has an~invariant interval domain~$U$, then the interior of~$\Sc\setminus U$ is~a~backward-invariant interval-domain, hence one can study separately the restriction of~this action to~$U$ and the restriction of~the inverse action to $\Sc\setminus U$. In the smooth case, one can apply results from \cite{KV} to both these restrictions, and prove that the action is somewhat regular.

The previous study of smooth semigroup actions on the circle has suggested the following paradigm: a \emph{generic} action either possesses a~forward-invariant interval-domain, or is two-sided minimal. In the latter case, even though there is an attraction, the attracting and repelling sections of the corresponding skew product are everywhere dense, see \cite{KN,A}.
Our results provide a~description for possible minimal sets of circle dynamics and establish, in one of the possible forms, the above paradigm. Namely, our first result is the following

\begin{theorem}
    \label{t:mins}
    Consider an iterated function system on the circle generated by orientation-preserving homeomorphisms $f_1,\dots,f_s$. Then
    \begin{itemize}
        \item either this system admits a forward-invariant interval-domain (and thus also a backward-invariant one);
        \item or there exist a unique forward-minimal set $\mM_+$, a~unique backward-minimal set $\mM_-$, and these sets have non-empty intersection.
    \end{itemize}
    Moreover, in the latter case there is an \emph{essential} intersection point of $\mM_+$ and $\mM_-$, in the following sense.
\end{theorem}
\begin{definition}
    An intersection point $x$ of the minimal sets $\mM_+$ and $\mM_-$ is called \emph{essential} if either both these sets accumulate to $x$ from the left, or both these sets accumulate to $x$ from the right.
\end{definition}

To state the second result, we need some definitions.
\begin{definition}
An interval-domain $U$ is \emph{strictly absorbing} for an IFS $F=(f_1, \dots, f_s)$ if for any its generator $f_i$ we have $f_i(U)\Subset U$.
\end{definition}
\noindent Note that the set of IFS's with a strictly absorbing interval-domain is open, because such domain survives under $C^0$-small perturbations. Let us denote the $C^r$-interior of the complement to this set by $\nd^r$. More precisely, by $\nd^r$ we denote the \emph{$C^r$-interior} of the set of IFS's $F=(f_1, \dots, f_s)$ such that
\begin{itemize}
\item $f_1, \dots, f_s$ are homeomorphisms (for $r=0$) or $C^r$-diffeomorphisms (for $r>0$) and preserve the orientation;

\item $F$ has no strictly absorbing interval-domain.
\end{itemize}

\noindent We shall say that an IFS is \emph{$C^r$-robustly minimal}, if it has a $C^r$-neighborhood consisting of minimal IFS.

Then the paradigm mentioned above is justified by the following

\begin{theorem}
    \label{t:C1}
    $C^1$-robustly forward and backward minimal IFS's are dense in $\nd^r$ for each $r = 1, \dots, \infty$.
\end{theorem}

\noindent Therefore, an alternative holds on an open and dense subset of the set of $C^r$-IFS's for $r \ge 1$: either an IFS admits a strictly absorbing interval domain, or it is robustly minimal.

Let us also add that by \autoref{c:diophantine} below IFS's such that some composition of $f_i$ is conjugated to an irrational rotation form a dense subset of $\nd^r$ for each $r \ge 0$.

Note that minimal systems have no internal points in the space of all $C^0$ IFS, as any system can be perturbed by making a Denjoy-like explosion, replacing points of a given orbit of a group, generated by this IFS, with (arbitrarily small) intervals (extending the maps to these intervals in an affine way). Though, Theorem~$\ref{t:C1}$ together with a simple observation from general topology will imply the following
\begin{proposition}
  \label{residual}
    Forward and backward minimal IFS's form a residual subset of $\nd^0$.
\end{proposition}

Now, let us introduce randomness to the dynamics.
\begin{definition}
A \emph{random dynamical system} (RDS), associated to an IFS $F=(f_1,\dots,f_s)$, is given by the probabilities $p_1,\dots,p_s$ that we associate to the corresponding maps. We assume that $p_i$ satisfy
$$
p_1+\dots+p_s=1, \quad p_1,\dots,p_s>0.
$$
\end{definition}
\todo[inline]{add a few words!}

\begin{definition}
A RDS is called \emph{globally synchronizing} if for any two points $x,y\in\Sc$ the distance between their random trajectories tends to~$0$ almost surely. Let us call an IFS \emph{globally synchronizing} \todo{added GS IFSs}if any RDS associated to this IFS is globally synchronizing.
\end{definition}

Strongly using Malicet's results (\cite{M}), we prove the following
\begin{theorem}
    \label{t:synch}
        Globally synchronizing IFS's form an an open and dense subset of $\nd^r$ for each $r = 0, \dots, \infty$.
\end{theorem}

We also discuss not globally synchronizing systems and the related property of \emph{factorization} (see \autoref{d:factor} below) in \autoref{s:ngss}.

We will conclude our paper (\autoref{s:examples}) with several examples illustrating different possibilities for the IFS's on the circle. The first of them is an example of a semi-group action that is minimal backwards, but not forward. Due to the existence of such example, whenever we use both forward and backward minimality, neither of this conditions can be omitted. Another examples show that even if there is no absorbing interval-domain, minimal sets can be different but both nowhere dense, or even intersecting only by a finite number of points (so that the conclusion of Theorem~\ref{t:mins} cannot be strengthened). The last example addresses the factorization property and necessary conditions for it to be claimed (see section~\ref{s:q-factor}). Finally, we show that in absence of an absorbing interval-domain, \emph{a priori}, one cannot claim factorization neither on~$\mM_-$, nor on~$\mM_+$.


%% file: preliminaries.tex
\section{Preliminaries}
\subsection{Notation}
Consider an IFS $F=(f_1, \dots, f_s)$ on the circle.
We will always assume that the generators $f_i$ are orientation preserving circle homeomorphisms.

Denote by $Σ^s$ the space of \emph{one-sided infinite sequences} $ω=ω_0ω_1\dots$ with $ω_i\in\set{1, 2, \dots, s}$, and by $W^s$ the set of all \emph{finite words} $w = w_0 \dots w_k$, $w_i\in\set{1, 2, \dots, s}$.
We shall use the following notation
\begin{align*}
    f_{ω, n}&≔f_{ω_{n-1}}\circ\dots\circ f_{ω_0},&ω&∈Σ^s;\\
    f_{w, n}&≔f_{w_{n-1}} \circ \dots \circ f_{w_0},&w&∈W^s,\,|w|≥n;\\
    f_w&≔f_{w,|w|},&w&∈W^s.
\end{align*}
Here $|w|$ is the length of a finite word $w$.
By $\Gr(F)$ we denote the semigroup generated by $f_1, \dots, f_s$.

\subsection{A measure in the base and stationary measures in the fiber}
Though the results we prove are purely topological, we use a Bernoulli probability measure to prove them.
From now on, we fix some positive numbers $p_i>0$, $p_1+\dots+p_s=1$.
The particular choice of $p_i$ will not be important for us.
Recall that the corresponding Bernoulli measure $\Pbase$ on $Σ^s$ is given by
\[
    \Pbase(ω_0=w_0, \; …, \; ω_k=w_k)=p_{w_0}…p_{w_k}
\]
on cylinders $ω_0=w_0,…,ω_k=w_k$.

Recall that a measure $μ$ is called \emph{invariant} for a map $f:S^1→S^1$ if $f_*μ=μ$, where $f_*μ(I)=μ(f^{-1}(I))$.
For an IFS, the fiber maps $f_i$ usually have no common invariant measure, so it is natural to consider a measure which is equal to the mean of its images.
\begin{definition}
    Let $F$, $p_i$ and $\Pbase$ be as above.
    A measure $μ$ is called \emph{forward stationary} for the random dynamical system $(F, \Pbase)$ if
    \[
         μ= p_1 \times (f_1)_*μ \;+\; \dots \;+\; p_s \times (f_s)_*μ.
    \]
    We say that a measure $μ$ is \emph{backwards stationary} for $(F, \Pbase)$ if it is forward stationary for the IFS $(f_1^{-1}, \dots, f_s^{-1})$ with the same distribution of probabilities.
\end{definition}

We recall that stationary measures describe the random orbit of any initial point.
\begin{lemma}
    [{see~\cite[Proposition~$3.17$]{M}}]
    \label{l:tostat}
    For any $x_0 \in \Sc$ and $\Pbase$-a.e. $\omega \in \Sigma^s$, the set of weak-$*$ accumulation points of the sequence of probability measures
    \[
        \mu_N=\frac 1 N \sum_{n=0}^{N-1} \delta_{f_{ω,n}(x_0)}
    \]
    consists of stationary probability measures.
\end{lemma}

\subsection{Factorization} \label{s:prel-factor}
In this section we discuss the following theorem due to D.~Malicet.
\begin{theorem}
    [{\cite[Theorem D]{M}}]
    \label{t:M}
	Let $F$ be an IFS such that its generators $f_1, \dots, f_s$ have no common invariant measure. Then (for any RDS associated to $F$) \todo{added ``for any RDS''}there exist a number $d$ and a~tuple of measurable functions $r_i\colon \Sigma^s \to \Sc$, $i=1, \dots, d$, such that for almost any $\omega \in \Sigma^s$ any closed interval $I$ included in $\Sc \setminus \set{r_1, \dots, r_d}$ is contracted by the corresponding random iterations:
	$\diam(f_{\omega, n}(I)) \to 0$ as $n \to +\infty$.
\end{theorem}

Note that for $d=1$ the system is globally synchronizing. A trivial example of IFS's with $d>1$ is given by \emph{factorizable} systems:
\begin{definition} \label{d:factor}
    We say that an IFS $F$ is \emph{factorizable} if there exist a number $d>1$ and a homeomorphism $R:S^1→S^1$ conjugated with the rotation $x↦x+\frac 1d$ such that $R$ commutes with $F$.
    
We say that an $F$ is \emph{factorizable on} a forward-invariant set $M$, if there exist a number $d>1$ and a homeomorphism $R:S^1→S^1$ conjugated with the rotation $x↦x+\frac 1d$ such that $R(M) = M$ and $R|_M$ commutes with $F|_M$.
\end{definition}
\noindent Clearly, a factorizable system cannot be globally synchronizing. See also the discussion in \autoref{s:questions} about the connection between factorization and $d>1$.
\todo[author=YK,inline]{Add more text explaining why it is important.}

\begin{remark} \label{r:d}
\todo{new}
For each $\omega \in \Sigma^s$ we can set as $d(\omega)$ the minimal number such that there exist $d(\omega)$ points $r_i \in \Sc$ such that any closed interval $I$ included in $\Sc \setminus \set{r_1, \dots, r_d}$ is contracted. Since this function is invariant, by the ergodic theorem it is constant  $\Pbase$-a.e. We will assume that $d$ is equal to this constant.   
\end{remark}
\begin{remark} \label{r:CIM}
    The set of IFS's with no common invariant measure is open in the space of all $C^r$ IFS's. Indeed, if $F_n\to F$ and $\mu_n$ is a common invariant measure of $F_n$, then any weak-$*$ limit point of $\set{\mu_n}$ is a common invariant measure for $F$. Further, this set is everywhere dense, since in any neighborhood we can find an IFS $F$ such that $f_1$ and $f_2$ are smooth Morse-Smale maps with disjoint sets of fixed points. Also note that under this condition any forward- or backwards-stationary measure has no atoms. Indeed, otherwise the set of atoms of maximal measure would be a common invariant finite set.
\end{remark}

%% file: questions.tex
\section{Questions} \label{s:questions}
There are several interesting questions, seeming to be really within the reach, that are related to the results of~\cite{M} and of this paper. 
The first of them concerns the factorization. Namely, the example presented in~\autoref{s:q-factor} below shows that in absence of a common invariant measure and absorbing interval domains, having $d>1$ (recall that $d$ was introduced in~\autoref{t:M}) does not guarantee factorization on $\mM_+$ (and hence, taking the inverse system, nor on $\mM_-$).

\begin{question}
    What is the best factorization-like statement that holds?
\end{question}

Note, that the same $d>1$ construction can be applied to the example with $\mM_+\cap \mM_-=\{0,\infty\}$. The same adding of identity to the semigroup and perturbing it only away from $\mM_-$ or only away from $\mM_+$ allows to create other non-factorizable examples with $d>1$.

\begin{question}
    What happens if we perturb the dynamics on both $\mM_+$ and $\mM_-$? What is the “maximal” perturbation that we can create without destroying $d>1$?
\end{question}

In~\cite{GKKV}, it was shown that for two-sided minimal dynamics, translation numbers of all the compositions in the semigroup determine the action up to the best possible extent: up to a conjugacy. This generalizes an analogous statement for groups, obtained by Ghys~\cite{G}. In this result, the action is assumed to have no invariant measure, and the conclusion is that it can be reconstructed up to a semi-conjugacy; it is also the best possible conclusion in this setting, as making a Denjoy-like explosion of an orbit does not change any rotation number.

\begin{question}
    \label{q:semi}
    What can be said in a semi-group setting without any minimality assumption?
\end{question}

For instance, can we distinguish between presence and absence of absorbing interval-domains? A uniform upper bound on the denominator does not suffice: for the canonical action of $\Gamma_2$, an index 6 subgroup of $\PSL(2,\bbZ)$, on the circle any map has zero rotation number, and hence an integer translation one. Though, three very strong NS maps, rotated third of a turn w.r.t. each other, have an absorbing domain formed by three intervals, and the rotation number of any element of the semigroup vanishes. Though, the translation number is not additive even on the semigroup, and hence one can not choose their lifts to the real line in such a way that the translation number of any composition vanishes.

As in Ghys' theorem, answering Question~\ref{q:semi} we cannot expect anything stronger than determining the IFS up to a semi-conjugacy. But even it, as the examples cited above show, cannot be claimed neither on $\mM_+$ (at least, unless $\mM_+\subset \mM_-$?), nor on $\mM_-$ (at least, unless $\mM_-\subset \mM_+$). So: what is the best statement that can be claimed?

\todo[inline]{(remove for now) Proposition \{ADD! Think: before or after?\} says that, unless the IFS possesses a common invariant measure, there is its neighborhood, in which the number d \{NAME\} is uniformly bounded. It means, that in a one-parametric family of IFS', in which no system has a common invariant measure, the d's that appear are uniformly bounded.}

Another interesting open-ended problem is the following 
\begin{question}
    How does the bifurcation between the “minimal” (or “globally contracting”) and “absorbing interval-domain” regimes occur?
\end{question}
\noindent In the case of a smooth dynamics, analogously to the famous Herman's results about the set of parameter values corresponding to an irrational rotation numbers, one can ask for the “bifurcation locus” between the two types of the behavior:
\begin{question}
    What can be said about the set of parameters, that correspond to neither minimal, nor interval-domains type systems for a generic family of $C^2$-systems? Is it of measure zero? Does it consist of finitely many points?
\end{question}

%% file: topology.tex
\section[Minimal sets]{The structure of the minimal sets}
This section is devoted to the proof of~\autoref{t:mins}. Consider an IFS $F=(f_1, \dots, f_s)$ that admits no absorbing interval-domain. Since the complement to a backwards-absorbing domain is forward-absorbing, our IFS has no backwards-absorbing domain either.
Let us first prove that $F$ has unique forward-minimal and backward-minimal sets, and then that these sets have non-empty intersection.
\subsection{Uniqueness of minimal sets}
Due to the symmetry, it is suffices to prove that the \emph{forward}-minimal set is unique. Suppose that there are at least two different forward-minimal sets, $\mM_+^1$ and $\mM_+^2$. Clearly, $\mM_+^1\cap\mM_+^2=\varnothing$. Let $U$ be the union of all connected components of $\Sc\setminus\mM_+^1$ that have non-empty intersection with $\mM_+^2$. 

\vspace{5mm}

\begin{tikzpicture}
\begin{scope}[shift={(0, 0)}, scale=0.2]

\draw (0, 0) circle[radius=6];
\foreach \an in {0, -5, -10, 30, 35, 55, 60, -130, -135, -140}
{
\draw (\an:6-0.5) -- (\an:6);
}
\node [below left] at (45:6) {$\mM_+^2$};

\draw (0, 0) circle[radius=10];
\draw [line width = 3pt] (-80:10) arc [radius=10, start angle=-80, end angle= 70];
\draw [line width = 3pt] (-170:10) arc [radius=10, start angle=-170, end angle= -100];
\node [below] at (-90:10) {$\mM_+^1$};
\node [right] at (15:10) {$I$};
\foreach \an in {-100, -95, -85, -80, 70, 75, 85, 90, -170, -175, -180}
{
\draw (\an:10-0.5) -- (\an:10);
}

\draw [dashed] (0:6) -- (0:10);
\draw [fill] (0:6) circle[radius=6pt];
\node [above right] at (0:6) {$c$};
\end{scope}

\begin{scope}[shift={(3, 0)}, scale=0.7]
\begin{scope}[shift={(0, -1.5)}, scale=1]
\draw (0, 0) -- (10, 0);
\draw [line width = 3pt, -] (2, 0) -- (8, 0);
\node [below] at (7, 0) {$f_j(I)$};

\draw [fill] (3, 0+2pt) circle[radius=2pt];
\node [above] at (3, 0) {$f_j(c)$};

\foreach \xs in {0.5, 2, 5, 8, 9}
{
	\begin{scope}[shift={(\xs, 0), (2, 0)}, scale=1]
	\draw (0, 0) -- (0, 0.2);
	\draw (0.2, 0) -- (0.2, 0.2);
	\draw (-0.2, 0) -- (-0.2, 0.2);
	\end{scope}
}
\node [below] at (0.5, 0) {$\mM_+^1$};

\draw (2+0.2, -1) -- (5-0.2, -1); 
\draw [fill] (2+0.2, -1) circle[radius=1pt];
\draw [fill] (5-0.2, -1) circle[radius=1pt];
\node [below] at (3.5, -1) {$I'$};
\draw [dashed] (2+0.2, 0) -- (2+0.2, -1);
\draw [dashed] (5-0.2, 0) -- (5-0.2, -1);
\end{scope}

\begin{scope}[shift={(0, 1.5)}, scale=1]
\draw (0, 0) -- (10, 0);
\draw [line width = 3pt, -] (3, 0) -- (7, 0);
\node [below] at (6, 0) {$I$};

\foreach \xs in {1, 3-0.2, 7, 8}
{
	\begin{scope}[shift={(\xs, 0), (2, 0)}, scale=1]
	\draw (0, 0) -- (0, 0.2);
	\draw (0.2, 0) -- (0.2, 0.2);
	\end{scope}
}
\node [below] at (1+0.2, 0) {$\mM_+^1$};

\draw [fill] (5, 0+2pt) circle[radius=2pt];
\node [above] at (5, 0+1pt) {$c$};

\end{scope}

\draw (5, -0.8) to [out=135,in=-135] (5, 1.5-0.7);
\draw (5-0.2, -0.8) -- (5, -0.8) -- (5, -0.8+0.2);
\node [left] at (4.7, 0) {$f_j$};
\end{scope}
\end{tikzpicture}

\noindent Since the distance between $\mM_+^1$ and $\mM_+^2$ is positive, $U$ is an interval-domain. Let us prove that $U$ is a~backwards-absorbing domain.

Denote by $\mI_0$ the set of connected components of $U$, put $N=|\mI_0|$. Consider an interval $I=(a,b)\in\mI_0$. Fix a point $c\in(a, b)\cap\mM_+^2$. By definition of $U$ and by forward-invariance of $\mM_+^{1,2}$, for any $j$ we have $f_j(a), f_j(b) \in \mM_+^1$, $f_j(c)\in \mM_+^2$. Therefore, $f_j(I)$ includes the interval $I'\in\mI_0$ that contains $f_j(c)$. For any $j$ the image $f_j(U)$ consists of $N$ disjoint intervals, and each of these intervals includes at least one subinterval from~$\mI_0$. Hence, all the intervals from $\mI_0$ are covered, and we have $f_j(U)\supset U$.

The contradiction proves that the minimal sets are unique.

\subsection{Intersection of \texorpdfstring{$\mM_+$}{ℳ₊} and \texorpdfstring{$\mM_-$}{ℳ₋}}
\begin{tikzpicture}

\draw (0, 0) -- (10, 0);
\draw [fill] (0.5, 0) circle[radius=2pt];
\node [above] at (0.5, 0) {$x$};
\node [above] at (5.6, 0) {$\mM_+$};
\node [below] at (9.5, 0) {$\mM_-$};
\node [below] at (4.2, 0) {$U_1$};

\foreach \n in {1, 3, 5, 7}
{
\draw (1.3^\n, 0) -- (1.3^\n, 0.2);
\draw (1.05*1.3*1.3^\n, 0) -- (1.05*1.3*1.3^\n, 0.2);
\draw (1.3*1.3^\n, 0) -- (1.3*1.3^\n, 0.2);
\draw (0.95*1.3^\n, 0) -- (0.95*1.3^\n, 0.2);
\draw [line width = 3pt, -] (1.3^\n, 0) -- (1.3*1.3^\n, 0);
}

\foreach \x in {2.5, 4.5, 7}
{
\draw (\x, 0) -- (\x, -0.2);
\draw (\x+0.2, 0) -- (\x+0.2, -0.2);
}
\foreach \x in {1.4, 8, 9}
{
\draw (\x, 0) -- (\x, -0.2);
}

\end{tikzpicture}

\noindent Let $\mI_1$ be the set of the intervals of the complement to $\mM_+$ that intersect $\mM_-$. Put $U_1:=\bigcup_{I\in \mI_1} I$.

Let us show that $U_1$ is a backwards-absorbing domain. Fix $j$ and an interval $I\in\mI_1$. Due to the forward-invariance of $\mM_+$ and backward-invariance of $\mM_-$, the interval $f_j^{-1}(I)$ does not intersect $\mM_+$ and intersects $\mM_-$. Therefore, the interval $J$ of $\Sc\setminus\mM_+$ that includes $f_j^{-1}(I)$ belongs to $\mI_1$. Thus, we have $f_j^{-1}(I)\subset U_1$. Since $I$ and $j$ are arbitrary, the domain $U_1$ is backwards-absorbing.

Recall that our iterated function system has no absorbing interval-domains, hence $|\mI_1|=\infty$.
Let $x$ be an accumulation point of the set of endpoints of these intervals.
Then either each left semi-neighborhood, or each right semi-neighborhood of $x$ includes infinitely many intervals from $\mI_1$. Since the closure of each interval $I\in\mI_1$ intersects both $\mM_+$ and $\mM_-$ (by the definition of $\mI_1$), the point $x$ is an essential intersection point of $\mM_+$ and $\mM_-$.

This concludes the proof of \autoref{t:mins}.

%% file: perturbation.tex
\section[Rotation number]{Perturbation of rotation numbers}

Consider an IFS $F=(f_1, \dots f_s)$ that admits no strictly absorbing interval-domain. 
Fix some lifts $\tilde f_1, \dots, \tilde f_s \colon ℝ→ℝ$ of the generators of $F$.
Put
\[
\tilde F+\delta = (\tilde f_1+\delta,\dots,\tilde f_s+\delta).
\]
\begin{proposition}
    \label{p:rotation-number}
If $F$ has no strictly absorbing interval-domain, then for any $\eps>0$ there exists a~word~$v$ such that $(\tilde f+\eps)_v(0) > (\tilde f-\eps)_v(0)+1$. 
\end{proposition}

\noindent This means that the maps $(\tilde f-\eps)_v$ and $(\tilde f+\eps)_v$ have different translation numbers, hence for some $\delta \in (-\eps, \eps)$ the translation number of $(\tilde f+\delta)_v$ is irrational. Since we can take any $\eps$, irrational translation number can be created by an arbitrary small perturbation. This observation will play the key role in the proofs of Theorems~\ref{t:C1} and~\ref{t:synch}. We add some technical details required for these proofs and state it as \autoref{c:diophantine} in section~\ref{s:diophantine} below.

\begin{remark}
For the special case when the generators of the IFS are projectivizations of linear maps from $SL(2, \mathbb R)$, a similar statement is known as \emph{Johnson’s Theorem} (\cite{J}).
\end{remark}

\begin{proof}[Proof of \autoref{p:rotation-number}]
Put $F^- = F-\eps, \; F^+ = F+\eps$ for shortness. For $\omega \in \Sigma^s$ put $I_{ω,n}=[\tilde f^-_{ω,n}(0), \tilde f^+_{ω,n}(0)]$.

Let $μ_-$ be any backward-stationary measure of $F^-$, let $\tilde μ_-$ be the (infinite) lift of $μ_-$ to the real line.

In Sections~\ref{s:dovorot} and~\ref{s:eps-min} below we will prove the following lemma.
\begin{lemma}
    \label{l:Exp}
    $\Exp \tilde \mu_-( I_{ω,n} )→∞$ as $n→∞$.
\end{lemma}
\noindent By this lemma for some $n$ and $ω$ we have $\tilde μ_-(I_{ω,n})>1$. Since $μ_-$ is a probability measure, this implies $|I_{ω,n}|>1$. Thus for the word $w=\omega_0 \dots \omega_{n-1}$ we have $\tilde f^+_w(0)>\tilde f^-_w(0)+1$. This proves \autoref{p:rotation-number} modulo \autoref{l:Exp}.
\end{proof}

\subsection{The growth in expectation} \label{s:dovorot}
\begin{proof}[Proof of \autoref{l:Exp}]
For an IFS $G$ and $n \le m$ by $g_{ω, n, m}$ we shall denote the composition 
$g_{ω_{m-1}} \circ \dots \circ g_{ω_{n}}$.
    Let us compare $\Exp \tilde μ_-(I_{ω,n})$ with $\Exp \tilde μ_-(I_{ω,n+k})$.
    To obtain $I_{ω,n+k}$, we iterate $I_{ω,n}$ by $\tilde f^-_{ω, n, n+k}$ and add the following extra interval:
    \[
        J_{ω,n,n+k}=I_{ω,n+k} ∖ \tilde f^-_{ω, n, n+k}(I_{ω,n})=
        \left(\tilde f^-_{ω, n, n+k}\left(\tilde f^+_{ω,n}(0)\right), \tilde f^+_{ω, n, n+k}\left(\tilde f^+_{ω,n}(0)\right)\right].
    \]
    Recall that $μ_-$ is a stationary measure for $F^-$, hence
    \[
        \Exp\tilde μ_-\left( \tilde f^-_{ω, n, n+k}(I_{ω,n}) \right) = \Exp\tilde μ_-(I_{ω,n}).
    \]
    Therefore,
    \[
        \Exp\tilde μ_-(I_{ω,n+k})-\Exp\tilde μ_-(I_{ω,n})=\Exp\tilde μ_-(J_{ω,n,n+k}).
    \]
    In particular, $\Exp\tilde μ_-(I_{ω,n+1})≥\Exp\tilde μ_-(I_{ω,n})$.

    Let us prove that for some $K∈ℕ$ and $δ>0$ we have
    \begin{equation}
        \label{eqn:Exp-mu-Ext}
        ∀n∈ℕ \;\; \Exp\tilde μ_-(J_{ω,n,n+K})≥δ.
    \end{equation}
    Obviously, this will imply \autoref{l:Exp}.
    Fix an interval $J⊂S^1$ such that $|J|=ε$ and $μ_-(J)>0$.
    We shall postpone the proof of the following purely topological lemma to the next subsection.

    \begin{lemma}
        \label{l:eps-covers}
        Let $F$, $ε$, $F^-$, $F^+$, $J$ be as above.
        Then for any $x∈S^1$ there exists a word $w$ such that $[f_w^-(x), f_w^+(x)]⊃J$.
        Moreover, there exists $K∈ℕ$ such that the word $w$ can be chosen to be at most $K$ symbols long.
    \end{lemma}

Let us prove \eqref{eqn:Exp-mu-Ext} for $K$ from \autoref{l:eps-covers} and $δ=\min p_i^K×μ_-(J)$. For $w \in W^s$ put $\Exp_v(\cdot) :=  \Exp \tilde μ_-(\cdot \; |ω_0…ω_{n-1}=v)$. It is enough to show that for any word $v$ of length $n$ we have
    \begin{equation}
        \label{eqn:Exp-mu-Ext-cond}
        \Exp_v \tilde μ_-(J_{ω,n,n+K})≥δ,
    \end{equation}   
    Fix $v$, then we know $I_n$.
    Taking as $x$ the right endpoint of $I_n$ and applying \autoref{l:eps-covers}, we get a word $w$ with $|w| \le K$. With probability at least $\min p_i^{|w|} \ge \min p_i^K$ we have $ω_n…ω_{n+|w|-1}=w$. In this case $J_{\omega, n, n+|w|}⊃J$, thus
\begin{equation}
	\label{eqn:exp-length-omega}
        \Exp_v \tilde μ_-(J_{ω,n,n+|w|})≥\min p_i^K×μ_-(J)=δ.
\end{equation}
Clearly, $J_{ω,n,n+K} \supset \tilde f^-_{ω, n+|w|, n+K}(J_{ω,n,n+|w|})$. 
Since $\mu_-$ is stationary, 
\[
	\Exp_v \tilde μ_-(J_{ω,n,n+|w|}) =
	\Exp_v \tilde μ_- \left( \tilde f^-_{ω, n+|w|, n+K}(J_{ω,n,n+|w|}) \right) \le
	\Exp_v \tilde μ_-(J_{ω,n,n+K}).
\]
    Together with~\eqref{eqn:exp-length-omega} this implies \eqref{eqn:Exp-mu-Ext-cond}, hence \eqref{eqn:Exp-mu-Ext}, and \autoref{l:Exp}.
\end{proof}

\subsection{Minimality by \texorpdfstring{$ε$}{ε}-orbits} \label{s:eps-min}
\begin{definition}
    For an IFS $G$, a sequence $x_0, \dots, x_k$ is called its \emph{$ε$-orbit}, if
    \[
        x_1 = g_{i_0}(x_0) + \eps_0, \dots, x_k = g_{i_{k-1}}(x_{k-1}) + \eps_{k-1},
    \]
    where $|\eps_j| < ε$.
\end{definition}
\noindent Recall that the IFS $F$ from \autoref{l:eps-covers} has no strictly absorbing interval-domain.
\begin{lemma}
    \label{l:eps-transitivity}
    Suppose that an IFS has no strictly absorbing interval-domain.
    Then for each $ε>0$, $x, y \in \Sc$ there is an $ε$-orbit that starts at $x$ and ends at $y$.
\end{lemma}
\begin{proof}
    Fix $x$ and $ε$.
    Let $U$ be the set of points $y$ such that there is an $ε$-orbit from $x$ to $y$.
    It is easy to see that $U$ is an open set, each connected component of $U$ has length at least $2ε$, and for each fiber map $f_i$, $U$ includes the $ε$-neighbourhood of $f_i(U)$.
    Thus either $U=S^1$, or $U$ is a strictly absorbing interval-domain.
    This immediately implies the lemma.
\end{proof}

\begin{proof}
    [Proof of \autoref{l:eps-covers}]  
    By \autoref{l:eps-transitivity} we may find an $(ε/2)$-orbit of $F$ that goes from $x$ to the center of $J$:
    \[
        x_0=x, \; x_1 = f_{w_0}(x_0) + \eps_0, \dots, x_k = f_{w_{k-1}}(x_{k-1}) + \eps_{k-1},
    \]
    \[
    J = (x_k-\eps/2, x_k+\eps/2).
    \]
    Moreover, by a compactness argument the number $k$ can be taken to be less than some constant $K$ independent of the initial point $x$. 
    Consider the word $w = w_0 \dots w_{k-1}$.
    Clearly, $x_{k-1}∈[f^-_{w,k-1}(x), f^+_{w,k-1}(x)]$, hence
    \begin{align*}
        f^-_w(x)&≤f_{w_{k-1}}(x_{k-1})-ε≤x_k-\frac{ε}{2},&
        f^+_w(x)&≥f_{w_{k-1}}(x_{k-1})+ε≥x_k+\frac{ε}{2},
    \end{align*}
    thus $[f_w^-(x), f_w^+(x)]⊃J$.
\end{proof}

\subsection{Obtaining irrational rotations} \label{s:diophantine}

\begin{corollary}
    \label{c:diophantine}
    For each $r = 0, \dots, \infty$ there is a dense subset $\rot^r \subset \nd^r$ such that for any IFS $G=(g_1, \dots, g_s) \in \rot^r$
    \begin{itemize}
        \item $G$ is a $C^\infty$-IFS;
        \item all generators of $G$ are Morse-Smale;
        \item there is a word $w \in W^s$ such that the map $g_w$ is $C^\infty$-conjugated to a pure irrational rotation. 
    \end{itemize}
    \emph{Moreover}, for any word $u \in W^s$ there is a dense subset $\rot^r_u \subset \nd^r$ of IFS's as above such that the word $u$ is a prefix of $w$.
\end{corollary}

\begin{proof}
    For a word $u$ and an open set $𝒰⊂\nd^r$, let us find an IFS $G ∈ 𝒰 ∩ \rot^r_u$.
    Fix a $C^\infty$ IFS $F∈𝒰$ such that all its generators are Morse-Smale.
    Take $ε>0$ such that for any $\delta \in [-ε, ε]$ the IFS $F+\delta$ still has Morse-Smale generators and belongs to $𝒰$.
    Then $F+δ$ satisfies the first two requirements for all $δ∈[-ε, ε]$.

    Let us prove that $(f+\delta)_w$ is conjugated to an irrational rotation for some $\delta \in [-ε, ε]$ and some word $w$ starting with $u$.
    As above, let $\tilde F$ be a lift of $F$ to the real line.
    By \autoref{p:rotation-number}, there is a word $v$ such that $\tilde f^+_v(0) > \tilde f^-_v(0)+1$.
    Take $w=uv$ and $p=\tilde f_u^{-1}(0)$. Clearly, $\tilde f^+_u(p) > 0 > \tilde f^-_u(p)$. Thus
\[
    \tilde f^+_w(p)  = \tilde f^+_v(\tilde f^+_u(p)) > \tilde f^+_v(0) > \tilde f^-_v(0) + 1 > \tilde f^-_v(\tilde f^-_u(p)) + 1 = \tilde f^-_w(p) -1,
\]    
    so $\tilde f^-_w$ and $\tilde f^+_w$ have different translation numbers.
    Hence, for some $\delta \in (-ε, ε)$ the translation number of $(\tilde f+\delta)_w$ is diophantine, and therefore due to the Herman's result~\cite{H} $(f+\delta)_w$ is $C^\infty$-conjugated to a rotation.
    Thus $F+\delta\in \mathcal U \cap \rot^r_u$.
\end{proof}

%% file: minimality.tex
\section{Minimality} \label{s:min}
The following well-known mechanism to prove the robust minimality was used, e.g., by Gorodetski and Ilyashenko~\cite{GI}.
\begin{proposition}
    \label{p:Hut}
    The semigroup generated by an irrational rotation~$R$ and $f\in C^1(\Sc)$ with a hyperbolic attracting fixed point~$a$ is $C^1$-robustly forward minimal.
\end{proposition}
\noindent
For reader's convenience, we provide the proof of this proposition in \autoref{s:robust}.

\begin{proof}
    [Proof of \autoref{t:C1}]
    Take the dense subset $\rot^r \subset \nd^r$ from \autoref{c:diophantine}.
    This set consists of robustly minimal IFS's by \autoref{p:Hut}.
\end{proof}


\begin{proof}
    [Proof of \autoref{residual}]
    For any open arc $I \subset \Sc$ we consider the set $\nd_I \subset \nd^0$ of IFS's such that forward orbit of any point intersects $I$.
    This set is open, since for such IFS circle is covered by finite number (by compactness) of preimages of $I$, and this property is robust.
    Note that $\nd_I$ is dense in $\nd^1$.
    Indeed, by \autoref{t:C1} minimal IFS's are dense in $\nd^1$, and any minimal IFS belongs to $\nd_I$. Since $\nd^1$ is dense in $\nd^0$, the set $\nd_I$ is dense in $\nd^0$.
    Intersecting the sets $\nd_I$ over all arcs $I$ with rational endpoints, we get the desired residual set.
\end{proof}

%% file: factorization.tex
\section[Not synchronizing systems]{Not globally synchronizing systems} \label{s:ngss}

In this section we will only consider IFS's without an absorbing domain such that its generators have no common invariant measure. We are interested in the case when the number $d$ from \autoref{t:M} is greater than one. By \autoref{t:mins} such IFS has a unique forward-invariant set $\mM_+$ and a unique backward-invariant set $\mM_-$. In \autoref{p:factor} below we shall prove that if $\mM_- = \Sc$, the IFS is factorizable (in the sense of \autoref{d:factor}) on $\mM_+$, and vice versa. But first we need to state some preliminary facts related to \autoref{t:M}.

Proposition~4.8 in~\cite{M} claims that supports of different ergodic stationary measures are disjoint and coincide with the forward-minimal compacts.
So in our situation there exists a unique forward-stationary probability measure $\mu_+$ with $\supp \mu_+ = \mM_+$. Similarly, there exists a unique backward-stationary probability measure $\mu_-$ with $\supp \mu_- = \mM_-$.

\autoref{t:M} provides us with $d$ functions $r_i: \Sigma^s \to \Sc$ (defined $\Prob$-a.e.) such that any circle arc in the fiber over $\omega$ that is disjoint from $\cup r_i(\omega)$ is contracted. For $\omega \in \Sigma^s$ set $R(\omega) := \set{r_1(\omega), \dots, r_d(\omega)}$.

\begin{lemma} \label{l:factor}
For $\Prob$-a.e. $\omega \in \Sigma^s$,

\begin{enumerate}
\item \label{a:1d}
The points $r_1(\omega), \dots, r_d(\omega)$ divide the circle into $d$ arcs with $\mu_+$-measure $1/d$ each,

\item \label{a:inv}
$R(\sigma(\omega))=f_{\omega_0}(R(\omega))$,

\item  \label{a:r-a-mu}
$\frac 1 d ((r_1)_*\Prob + \dots + (r_d)_*\Prob) = \mu_-$.
\end{enumerate}
\end{lemma}

\begin{proof}
Conclusion~\ref{a:1d} is proven in~\cite[Proposition~$4.3$]{M}.

Since $f_{\omega, n+1} = f_{\sigma(\omega), n} \circ f_{\omega_0}$, an interval $I$ is contracted by the sequence $f_{\omega, k}$ iff the interval $f_{\omega_0}(I)$ is contracted by the sequence $f_{\sigma(\omega), k}$. Since any interval that lies in $\Sc \setminus R(\omega)$ is contracted by $f_{\omega, k}$, any interval in $\Sc \setminus f_{\omega_0}(R(\omega))$ is contracted by $f_{\sigma(\omega), k}$. Hence, by~\autoref{r:d} we have $R(\sigma(\omega))=f_{\omega_0}(R(\omega))$.

Using conclusion~\ref{a:inv}, one may check that the measure $\frac 1 d ((r_1)_*\Prob + \dots + (r_d)_*\Prob)$ is backward-stationary. Hence it coincides with the unique backward stationary probability measure $\mu_-$.
\end{proof}

\begin{proposition} \label{p:factor}
Let $G=(g_1, \dots, g_s)$ be a backward-minimal IFS such that $d>1$. Then the IFS is factorizable on $\mM_+$.
\end{proposition}
\begin{proof}
For convenience, we will prove the proposition for the inverse IFS $F = (f_1, \dots, f_s) = (g_1^{-1}, \dots, g_s^{-1})$. That is, let us prove that if the system is forward-minimal and $d>1$, then the inverse IFS is factorizable on $\mM_-$.

As discussed above, $F$ has a unique forward-stationary measure $\mu_+$ and a unique backward-stationary measure $\mu_-$. By \autoref{r:CIM} the measure $\mu_+$ has no atoms.
Since the IFS is forward-minimal, $\supp \mu_+ = \Sc$.
Hence, one can define the $1/d$-rotation $T$ with respect to $μ_+$, i.e. a map such that $μ_+([x, T(x)])=\frac 1 d$ for all $x∈S^1$.

In order to show that $F$ is factorizable (with respect to the homeomorphism $T$), we will apply~\autoref{l:factor}. By conclusion~\ref{a:1d} (of \autoref{l:factor}) the map $T$ permutes the terms of the formula from conclusion~\ref{a:r-a-mu}. Hence $T$ preserves $\mu_-$ and, therefore, $T$ also preserves $\mM_- = \supp \mu_-$.

By conclusion~\ref{a:inv}, $\Prob$-a.e. we have
\begin{align*}
R(\sigma(\omega)) &= f_{\omega_0}(R(\omega)), \\
f_{\omega_0}^{-1}(R(\sigma(\omega))) &= R(\omega).
\end{align*}
By conclusion~\ref{a:1d} $T$ cyclically permutes $\{ r_i \}$, so for every $i$
\[
T \circ f_{\omega_0}^{-1}(r_i(\sigma(\omega))) = f_{\omega_0}^{-1} \circ T(r_i(\sigma(\omega)))    \qquad \Prob-\text{a.e.}.
\]
Clearly, $\omega_0$ and $\sigma(\omega)$ are independent and $\sigma_*\Prob = \Prob$. So for each $j$ and $\Prob$-a.e. $\tilde \omega$ we have
\begin{align}
T \circ f_j^{-1}(r_i(\tilde \omega)) &= f_j^{-1} \circ T(r_i(\tilde \omega)), \nonumber \\
\label{e:commute}
T \circ f_j^{-1}(x) &= f_j^{-1} \circ T(x) \qquad \qquad (r_i)_*\Prob \text{- a.e.}
\end{align}
Taking average over $i$, by  conclusion~\ref{a:r-a-mu} we get that~\eqref{e:commute} holds $\mu_-$-a.e. Hence~\eqref{e:commute} holds on a dense subset of $\mM_-$ and, by continuity, on the whole $\mM_-$. This finishes the proof of \autoref{p:factor}.

\end{proof}

%% file: synchronization.tex
\section{Synchronization}\label{sec:sync}
\subsection{From robust non-synchronization to robust factorization}
In this section we prove~\autoref{t:synch}. \cite[Corollary~$2.9$]{M} states that global synchronization is a $C^0$-open property. Let us prove that this property is dense in $\nd^r$. Since $C^\infty$-IFS's are dense in $\nd^r$ for any $r$, we may and will assume that $r=\infty$.

 We argue by a contradiction. Suppose that there is an open set $U_0 \subset \nd^\infty$ consisting of non-synchronizing IFS's. Let us find an open $U_1 \subset U_0$ such that any IFS in $U_1$ is factorizable for some $d>1$. To do so, we take as $U_1$ any open subset of $U_0$ such that for any IFS in $U_1$
\begin{itemize}
\item its generators have no common invariant measure;
\item it is forward and backwards minimal.
\end{itemize} 
Such set exists, because both conditions are open and dense in $\nd^\infty$: the first by~\autoref{r:CIM}, and the second by~\autoref{t:C1}. 

For any IFS $G \in U_1$ let us prove that it is factorizable. Consider the number $d$ from~\autoref{t:M} for the IFS $G$. If $d=1$, the IFS $G$ would be globally synchronizing. But $G \in U_0$, and $U_0$ consists of non-synchronizing IFS's. Thus $d>1$. \autoref{p:factor} claims factorization on the unique forward-minimal set $\mM_+$. Since $G$ is forward minimal, $\mM_+ = \Sc$, and $G$ is factorizable. 

\subsection{There is no open set of factorizable IFS's}
To finish the proof of \autoref{t:synch}, we need the following two lemmas.
\begin{lemma} \label{l:smooth-factor}
    Let $F$ be a $C^∞$-smooth factorizable IFS, $R\circ f=f\circ R$ for $f∈\Gr(F)$.
    Suppose that $\Gr(F)$ contains a map smoothly conjugated to an irrational rotation.
    Then $R$ is a smooth map.
\end{lemma}
\begin{proof}
    Let $f∈\Gr(F)$ be a map smoothly conjugated to an irrational rotation.
    Let us write all maps in the chart such that $f(x)=x+α$.
    Then we have $R(x+α)=R(x)+α$, hence $R(nα)=R(0)+nα$, $n∈\bbZ$.
    Since $α$ is irrational, the set $\set{nα|n∈\bbZ}$ is dense in $S^1$, hence $R$ is the rotation $R(x)=R(0)+x$.
    Finally, $R$ is smooth in the linearizing chart for $f$, hence it is smooth in the original chart as well.
\end{proof}

\begin{lemma}
    \label{l:mult}
    There is an open and dense subset $A \subset \nd^∞$ such that for any IFS $G \in A$ there is a map $g∈\Gr(G)$ such that
    \begin{itemize}
        \item $g$ is a Morse-Smale map with zero rotation number;
        \item all attractors of $g$ have different multipliers.
    \end{itemize}
\end{lemma}
\noindent This lemma looks rather obvious, but its proof is a bit technical. So we separate this proof in Appendix~B. 

Let us continue the proof of \autoref{t:synch}. By \autoref{l:mult} there are an open $U_2 \subset U_1$ and a word $w$ such that for any IFS $G \in U_2$ the map $g_w$ is a Morse-Smale map with zero rotation number, and the multipliers at all its fixed points are different. By~\autoref{c:diophantine} we may find an IFS $G \in U_2$ such that $\Gr(G)$ contains a map that is smoothly conjugated to an irrational rotation. Take any attractor $a$ of the map $g_w$. Since $G$ is factorizable, $g_w\circ R=R\circ g_w$ for some $R$. Thus the point $R(a)$ is another attractor of $g_w$. Since $R$ is smooth (by \autoref{l:smooth-factor}), it has the same multiplier as $a$. This contradiction proves \autoref{t:synch}.

%% file: examples.tex
\section{Zoo of examples} \label{s:examples}
\subsection{Minimality in one direction, not in the other}\label{s:one-sided-min}
Take a Schottky group $G=\langle f_1,f_2 \rangle$, corresponding to a torus with a hole (that is, a group, generated by a strong North--South map $f_1$ and a strong East--West map $f_2$). It has a Cantor limit set $\Lambda$, and any interval of its complement may be sent to any other; moreover, the map $f:=[f_1,f_2]$ is a hyperbolic transformation, fixing one of such intervals (the quotient by its action wraps $I_0$, which is homeomorphic to a line, back to a circle that is the border of the hole). Denote this interval by $I_0$, and let $v$ be a vectorfield on $I$ for which $f$ is its time-one map ($f$ is a hyperbolic Moebius transformation, so we can find one).

Now, add to this group a transformation $f_3$, which is identity outside $I_0$, and equals to the time \emph{minus} $\sqrt{2}$ flow of the vectorfield $v$ inside it. Then, for any point $x\in I_0$ its $\langle f, f_3\rangle_+$-backwards orbit is dense in $I_0$ (as we have time-one and time-minus $\sqrt{2}$ shifts along the vector field~$v$), and hence its $\langle G, f_3\rangle_+$-backwards orbit is dense in all the circle. Finally, add a map $f_4$ that sends $\Lambda$ strictly inside itself: $f_4(\Lambda)\subsetneq \Lambda$. For instance, one can take a map that is piecewise-Moebius, glued from the contracting pieces of $f_1$ and $f_2$, extending it in any continuous way in between.

Then $\Lambda$ is forward-invariant under all the constructed maps, and is $G$-minimal. Hence, for $G_1:=\langle f_1, f_1^{-1}, f_2, f_2^{-1}, f_3, f_4 \rangle_+$ we have $\Lambda$ as a forward-minimal set. Let us now show that any backwards $G_1$-orbit turns out to be dense.

Note first that such orbit of any $x\in \Sc$ cannot be contained in $\Lambda$, finding $x'\in G_1^{-1}(x)\setminus \Lambda$.
Indeed, if $x\notin \Lambda$, we are already done; otherwise, take $y\in \Lambda\setminus f_4(\Lambda)$. Then, there exists a neighborhood $U\ni y$ such that $U\cap f_4(\Lambda)=\emptyset$, and hence that $f_4^{-1}(U)\subset \Sc\setminus \Lambda$. By minimality of the action of $G$ on $\Lambda$, there is $g\in G$ such that $g(x)\in U$, and hence $(f_4^{-1}\circ g)(x)\notin \Lambda$. Note that as we have taken all the four maps $f_1, f_1^{-1}, f_2, f_2^{-1}$ among the generators of $G_1$, the map $g$ belongs also to $G_1^{-1}$.

Finally, $G$-action (both backward and forward) allows to send any interval of the complement to $\Lambda$ to $I_0$. Hence, $G_1^{-1}(x)$ intersects $I_0$, and by the above arguments is dense there. Hence, the orbit $G_1^{-1}(x)$ is dense in any interval of complement to $\Lambda$, and hence in all the circle.

We have $\mM_+=\Lambda\subsetneq \mM_-=\Sc$.

\subsection{No minimality in either way, one of the minimal sets is a subset of the other}

Take the previous example. Choose a full (both forward and backward iterations authorized) orbit of a point $x_0\in \Sc$ and make a Denjoy-like construction. Namely, replace each point $g(x)$ of this orbit by a small interval $I_{g(x)}$ (ensuring that the sum of their lengths converges), and extend maps to these intervals in an affine way: $h\colon I_{g(x)}\to I_{h(g(x))}$. Take the constructed maps to generate a semi-group.

Consider the complement $K=\Sc \setminus \bigcup_{g} I_{g(x)}$ to the union of these intervals, together with the natural “gluing back” projection $\pi\colon\Sc \to \Sc$ that collapses the glued intervals and hence semi-conjugates the new dynamics to the original one. Then, the new minimal sets are
$$
\tM_{\pm}=\pi^{-1}(\mM_{\pm})\cap K,
$$
and hence
$$
\tM_{+}\subsetneq \tM_{-} \subsetneq \Sc.
$$

\subsection{No minimality in either way, the minimal sets are not subsets of each other: Cantorvals}

Consider a Denjoy example $f$ obtained by blowing up \emph{two} orbits of an irrational rotation. Then $f$ is minimal on a Cantor set $K$, and there are two families of disjoint intervals $I_i, J_i$ with $f(I_i)=I_{i+1}$, $f(J_i)=J_{i+1}$. As in the previous example, take a vector field $v$, nonzero on $I_0$ and on $J_0$ and vanishing outside them, and take two maps $f_1$ and $f_2$ to be time one and minus $\sqrt{2}$ maps along this vector field. Now, choose three intervals $I_{i'}=[p_1,p_2]$, $I_{i''}=[p_1',p_2']$, $J_k=[q_1,q_2]$ that are in the cyclic order $I_{i'},J_k,I_{i''}$. Finally, choose a homeomorphism $g$ so that
$$
g(p_1)=p_1, \quad g(p_2')=p_2'; \quad g([p_2,q_1])\subset J_k
$$
(see the figure), and $g$ is identity outside the arc $A:=[p_1,p_2']$.
\begin{figure}[h]
\centering
\begin{tikzpicture}
    \path
        coordinate (A1) at (0, 0)
        coordinate (A2) at (1, 0)
        coordinate (A3) at (2, 0)
        coordinate (A4) at (4.5, 0)
        coordinate (A5) at (5.5, 0)
        coordinate (A6) at ($ (A2)+(A5) $)
        coordinate (A10) at (0, -1)
        coordinate (B1) at ($ (A1)+(A10) $)
        coordinate (B2) at ($ (A3)+(A10)+(0.25, 0) $)
        coordinate (B3) at ($ (A3)+(A10)+(0.75, 0) $)
        coordinate (B4) at ($ (A4)+(A10)-(0.75, 0) $)
        coordinate (B5) at ($ (A4)+(A10)-(0.25, 0) $)
        coordinate (B6) at ($ (A6)+(A10) $);

    \draw[dotted,thick] (A2) -- (A3) (A4) -- (A5) (B2) -- (B3) (B4) -- (B5);
    \draw[very thick] (A1) -- (A2) (A5) -- (A6) (B1) -- (B2) (B5) -- (B6);
    \draw[very thick,dashed] (A3) -- (A4) (B3) -- (B4);

    \foreach \i in {1,...,6}{
        \fill (A\i) circle (2pt);
        \fill (B\i) circle (2pt);
        \draw[->] ($ (A\i)-(0, 0.15) $) -- ($ (B\i)+(0, 0.15) $);
    }
\end{tikzpicture}
\end{figure}

Then, for two Cantorvals $C_1=K \cup \left(\bigcup_i I_i\right)$ and $C_2=K \cup \left(\bigcup_i J_i\right)$ one has
$$
g(C_2)\subset C_2, \quad g(C_1)\supset C_1.
$$
Indeed, outside $A$ the map $g$ is the identity; at the same time, for any $x\in A$ at least one of $x\in I_{i'}\cup I_{i''}$ and $g(x)\in J_k$ holds.

Now, consider the semigroup
$$G:=\langle f, f^{-1},f_1, f_2, g \rangle_+.$$
Both sets $C_1,C_2$ are $f$-invariant, hence the set $C_2$ is $G$-forward invariant, while $C_1$ is $G$-backwards invariant. Finally, in the same way as in the previous example, we see that these sets are in fact minimal. Namely, for any point $x\in \Sc$ there is its $f$-image that lies on the arc $(p_2,q_1)$, and will be sent by $g$ to a point in $J_k$. Thus, the forward orbit of $x$ intersects $J_0$, and hence any $J_i$. Thus, it is dense in $C_2$. The same argument shows that the backwards orbit of any point is dense in $C_1$.

We thus have $\mM_+=C_2$, $\mM_-=C_1$.

\subsection{No minimality in either way, the minimal sets are not subsets of each other: Cantor sets}

This example is quite analogous to the previous one. Let
$$
T_{-1}, T_{0}, T_{1}\colon[-1,1]\to [-1,1]
$$
be the 100-times contractions to the points $-1$, $0$, $1$ respectively, and let $C$ be the Cantor set associated to this IFS. Then,
$$
C= \{0\} \cup \bigcup_{i=0}^{\infty} T_0^i (C_1) \cup \bigcup_{i=0}^{\infty} T_0^i (C_{-1}),
$$
where $C_{\pm 1}:=T_{\pm 1}(C)$.

Now, take the circle to be the projective line $\Sc=\bbR\cup\set{\infty}$, and extend the map $T_0(x)=x/100$ to all of it. Also, let $Q(x)=x/10$ (thus, $Q^2=T_0$) and define the Cantor sets
$$
K_+:= \set{0,\infty} \cup \bigcup_{i\in\bbZ} T_0^i (C_1\cup C_{-1}),
$$
$$
K_-:=Q(K_+)=  \set{0,\infty} \cup \bigcup_{i\in\bbZ} T_0^i (Q(C_1)\cup Q(C_{-1})).
$$

Note that $K_-\cap K_+=\{0,\infty\}$, and that on the arc $(0,+\infty)$ parts of these sets go in the order
$$
\dots, T_0^{i}Q(C_1), T_0^i(C_1),T_0^{i-1}Q(C_1),T_0^{i-1}(C_1), \dots
$$

Take the map $g$ to be defined in the following way:
\begin{equation}\label{eq:g}
    g(x)=
    \begin{cases}
        x, & x<Q^{-1}(-1) \, \text{or} \, x>1,\\
        T_1(x) & x\in [-1,1], \\
        h T_{-1}^{-1} h^{-1} (x), & x\in Q^{-1}([-1,T_{-1}(1)],\\
        \text{any continuous interpolation} & \\
        \text{between $Q(1)$ and $T_1(-1)$}, & x\in [Q^{-1}(T_{-1}(1)), -1].
    \end{cases}
\end{equation}
where
\[
    h(x)=
    \begin{cases}
        Q(x), & x\ge 0 \\
        Q^{-1}(x) & x< 0.
    \end{cases}
\]
(see the figure).

\begin{figure}[h]
\centering\includegraphics{Maps-1}
\end{figure}

The first line of~\eqref{eq:g} ensures that $g$ is the identity transformation to the left of $Q^{-1} (-1)$ and to the right of $1$. The second, that $K_+\cap [Q^{-1} (-1),1]=C$ is sent to $T_1C$, and hence $g(K_+)\subset K_+$. Finally, the last one ensures that for $K_-\cap [Q^{-1} (-1),1] = h(C)$ one has $g^{-1}(h(C))\subset h(C)$, and hence $g(K_-) \supset K_-$.

Now, denote $\tilde{g}(x):=-g(-x)$ and take the semigroup $G_1:=\langle T_0, T_0^{-1}, g, \tilde{g} \rangle_+$. The previous arguments imply immediately that
$$
    G_1(K_+)\subset K_+, \quad G_1^{-1}(K_-)\supset K_-.
$$
Also, $g|_{[-1,1]}=T_1$, $\tilde{g}|_{[-1,1]}=T_{-1}$, and it is easy to see that $K_{\pm}$ are indeed the forward and backward minimal sets of $G_1$.

\subsection{Not so easy to factorize}\label{s:q-factor}

Take the example $G_1$ of a backward, but not forward minimal semigroup action on the circle. Lift it to a $d$-leaves cover, obtaining the action~$\tilde{G}_1$. Now, perturb it outside $\mM_+$; namely, we add to~$\tilde{G}_1$ the map $g$ that is identity everywhere, except for one interval $J$ of complement to $\tilde{\mM}_+$. On one hand, for the new system we still have the same $d$, as after collapsing the intervals of complement to $\tilde{\mM}_+$ the dynamics is $d$-factorizable. On the other hand, it is not globally factorizable due to the perturbation. Indeed, initial lift $\tilde{G}_1$ is factorizable in a unique way, as the repellers $r_i$ from \autoref{t:M} are present everywhere due to the backwards minimality. Though, $g$ cannot be pushed through this factorization, as for two preimages of the same interval of complement to $\mM_+$ on one (namely, on $J$) the map $g$ is not identity, while on the other it~is.

In particular, this means that the dynamics cannot be factorized on $\tilde{\mM}_-=\Sc$, and thus in \cite[Proposition 4.3]{M} one cannot replace forward-stationary measure with the backward stationary one. Also, reversing the maps we get an example of an action that cannot be factorized on $\mM_+$, and for which the forward-stationary measure $\mu$ of the attractor-attractor arc can be different from $1/d$.

%% file: robust.tex
\section{Robust minimality}\label{s:robust}
This section is dedicated to the proof of \autoref{p:Hut}. We use a well-known argument going back to Gorodetski and Ilyashenko,~\cite{GI}. The proof is based on the following classical result 
\begin{lemma*}[Hutchinson,~\cite{Hut}]
Consider two intervals $J \subset I$ and an IFS $F$. Suppose that there exists a~tuple of elements $g_i\in\Gr(F)$ such that 
    \begin{enumerate}
    \item each $g_i$ contracts in $I$,
    \item  $g_i(I) \subset I$ for each $g_i$,
    \item $J⊂⋃_i g_i(J)$.
    \end{enumerate}
    Then the action is minimal on~$J$, i.\,e., the forward $F$-orbits of each point $x\in J$ is dense in $J$.
\end{lemma*}
\begin{proof}
We shall prove that the IFS $G = \set{g_i}$ is minimal on $J$, this will imply that $F$ is also minimal. To do so, it is enough to prove that for any open $U \subset J$ there is a word $u$ such that $g_u(J) \subset U$. 

By the third condition for any point in $J$ its preimage under some $g_i$ is in $J$ again. Let us pick any point $x \in U$ and take its preimage many times, choosing $g_i$ as above. This gives us an arbitrary long word $u$ such that $g_u^{-1}(x) \in J$. If $u$ is long enough, $g_u$ strongly contracts on $I$ by the first two conditions, hence $g_u(J) \subset U$.
\end{proof}

\begin{proof}
    [Proof of \autoref{p:Hut}]
    Let $I\ni a$ be a~small open neighborhood such that $f$ contracts in~$I$ and $J \ni a$ be an open subset of $I$. Then we may find a tuple of numbers $(n_j)$ such that the tuple $\set{R^{n_j}∘f}$ satisfies the assumptions of the previous lemma.
    We shall also pick a tuple of numbers $(m_j)$ such that the images of $J$ under $R^{m_j}$ cover $\Sc$,
    \begin{equation}
        \label{eqn:R-cover}
        ⋃_jR^{m_j}(J)=\Sc=⋃_jR^{-m_j}(J).
    \end{equation}

    Consider any semigroup $\Gr_\eps=\Gr(R_\eps, f_\eps)$ sufficiently $C^1$-close to the original one.
    For $ε$ small enough, the tuple $\set{R_ε^{n_j}∘f_ε}$ still satisfies the assumptions of the previous lemma, and \eqref{eqn:R-cover} holds for $R_ε$ instead of $R$.
    Let us prove that such $\Gr_\eps$ acts minimally on the circle.
    Consider a point $x\in\Sc$ and an open interval $U\subset\Sc$.
    Due to \eqref{eqn:R-cover}, there exist $i$ and $j$ such that $R_ε^{m_i}(x)\in J$ and $R^{-m_j}(U)\cap J\neq\varnothing$.
    Finally, due to the previous lemma, the orbit of $R^{m_i}(x)$ visits $R^{-m_j}(U)\cap J$, hence the orbit of $x$ visits~$U$.
\end{proof}

%% file: multipliers.tex
\section{Different multipliers}
\begin{proof}[Idea of the proof of \autoref{l:mult}]
    Given an arbitrary open subset $U⊂\nd^∞$, we want to find an IFS $F∈U$ with the properties described above.
    Take an IFS $G\in U$ with Morse-Smale generators $g_i$.
    Choose $k$ such that $\rho(g_1^k)=0$.
    Then $g_1^k$ is a Morse-Smale map, but an attractor $a_i$ of $g_1^k$ has the same multiplier as $g_1(a_i)$.

    We want to slightly change $g_1^k$ so that the attractors of the new map will have different multipliers.
    For this, we shall take $G$ such that one of the maps $g_w∈\Gr(G)$ is conjugated to an irrational rotation, hence some of its powers $g_w^m$ is close to the identity.
    It turns out that after an appropriate perturbation of $G$ for $n$ large enough, the map $g_1^{kn}\circ g_w^m$ satisfies the required properties.
\end{proof}

\begin{proof}[Proof of \autoref{l:mult}]
    Take an IFS $G \in U \cap \rot^\infty_2$ guaranteed by \autoref{c:diophantine}.
    Recall that by the definition of $\rot^\infty_2$, all generators of $G$ are Morse-Smale, and there is a word $w=2\dots$, such that $g_w$ is smoothly conjugated to an irrational rotation.
    Choose $k$ such that the map $g_1^k$ has zero rotation number.
    Let $a_i$ denote the attractors of $g_1^k$.

    Fix $m$ such that $g_w^m$ is so close to $\Id$ that for any $n>0$ the map $g_{u(n)}$, $u(n)=w^m1^{kn}$,
    \begin{itemize}
        \item is Morse-Smale;
        \item has the same number of attractors as $g_1^k$;
        \item these attractors $b_i=b_i(n)$ are close to $a_i$, namely $b_i$ is in the basin of attraction of $a_i$ under the map $g_1^k$.
    \end{itemize}
    This property still holds if we replace the IFS $G$ by any IFS in some neighborhood of $G$, denote this neighborhood by $V \subset U$.
    Note that for any IFS in $V$ we have $b_i(n) \to a_i$ as $n \to \infty$.

Consider the word $v=w^m=2\dots$.
Denote by $O(x, v, G)$ the collection of points $g_{v,1}(x), \dots, g_{v,|v|}(x)$.
Let us replace the IFS $G$ by its small pertubation that still lies in $V$, to obtain the following property:
the union of the orbits $O(a_i, v, G)$ is disjoint from the set of attractors $a_i$.
To do so, consider the family $G^\eps: g_2^\eps = g_2+\eps, g_i^\eps = g_i$ for $i \ne 2$.
Note that $a_j^\eps =a_j$, while all points in $O(a_i^\eps, v, G^\eps)$ move monotonously in $ε$ since $v$ starts with $2$.
Thus for any $i$, $j$ the point $a_j$ belongs to $O(a_i^\eps, v, G^\eps)$ for a nowhere dense set of parameters $\eps$.
Therefore $G^\eps$ is as required for some for some $\eps$.

Recall that $b_i(n) \to a_i$ as $n \to \infty$.
Hence $O(b_i(n), v, G) \to O(a_i, v, G)$ as $n→∞$ as well.
Thus for large $n$ the union of the orbits $O(b_i, v, G)$ is disjoint from the set of attractors $b_j$.

This property means that we can perturb the multiplier of any $b_i$, leaving the multipliers of other $b_j$ unchanged.
To do so, we perturb $g_2$ in a small neighborhood of $b_i$ so that $g_2'(b_i)$ changes while $g_2(b_i)$ is preserved.

Let us check that such perturbation changes the multiplier of $b_i$ and preserves the multipliers of $b_j$ for $j≠i$.
For each $j$ (including $j=i$) the perturbation preserves the values of $g_l$ at all the points $b_j$ and at the union of the orbits $O(b_j, v, G)$, hence it preserves these orbits.
The multiplier of $b_j$ equals
\[
    \prod_{l=0}^{|u|} g'_{u_l}(g_{u,l}(b_j)).
\]
The first factor $g'_2(b_j)$ is preserved for $j≠i$ and changes for $j=i$.
The next $|v|-1$ factors correspond are derivatives of the fiber maps at $O(b_j, v, G)$, and thus are preserved by our perturbation.
For the remaining $kn$ terms, $u_l=1$, hence they are preserved since $g_1$ is not perturbed.

Finally, we can perturb each multiplier $g_{u(n)}'(b_i(n))$ independently, hence we can make all of them different.
\end{proof}